\documentclass[10pt,twoside]{article}
\usepackage{graphicx}
\usepackage{amsmath}
\usepackage{amsfonts}
\usepackage{amssymb}

\usepackage{Latex-document}

\title{\bf  Random Walks in Random Environments\vskip 6mm}

\author{Ofer Zeitouni\vspace*{-0.5cm}\thanks{Departments of Electrical Engineering and
of Mathematics, Technion, Haifa 32000, Israel. E-mail: zeitouni@ee.technion.ac.il}}

\date{\vspace{-8mm}}

\begin{document}
\newtheorem{theorem}[equation]{Theorem}
\newcommand{\PP}{{\mathbb P}}
\newcommand{\EE}{{\mathbb E}}
\newcommand{\ZZ}{{\mathbb Z}}
\newcommand{\NN}{{\mathbb N}}
\newcommand{\won}{{\boldsymbol 1}}
\newcommand{\calV}{{\mathcal V}}
\newcommand{\om}{\omega}
\newcommand{\Om}{\Omega}
\newcommand{\uu}{\underline}
\newcommand{\oo}{\overline}
\newcommand{\eps}{\epsilon}
\newcommand{\beps}{{\scriptscriptstyle{{ \cal E}}}}

\newcommand{\dsum}{\displaystyle\sum}
\newcommand{\dprod}{\displaystyle\prod}
\newcommand{\limn}{\lim_{n \rightarrow \infty}}

\maketitle

\thispagestyle{first} \setcounter{page}{117}

\begin{abstract}\vskip 3mm
Random walks in random environments (RWRE's)
have been a source of surprising
phenomena and challenging problems
since they began to be studied in the 70's.
Hitting times and, more recently, certain regeneration
structures, have played a major role in our understanding of RWRE's.
We review these and provide some hints on current research directions
and challenges.

\vskip 4.5mm


\end{abstract}

\vskip 12mm

\section{Introduction} \label{section 1}\setzero
\vskip-5mm \hspace{5mm }

Let $S$ denote the 2d-dimensional simplex,
set $\Omega=S^{Z^d}$,
and let
$\omega(z,\cdot)=\{\omega(z,z+e)\}_{e\in Z^d, |e|=1}$ denote the
coordinate of $\omega\in\Omega$ corresponding to $z\in Z^d$.
$\Omega$ is
an ``environment'' for an inhomogeneous  nearest neighbor
 random walk (RWRE) started at $x$ with
{\it quenched} transition probabilities
$P_\omega(X_{n+1}  = z+ e | X_n = z) = \omega
 (x,x+e)$ ($ e\in  {\ZZ}^d,
|e|=1$), whose law is denoted
$P^x_\om$. In the RWRE model, the environment is random, of law $P$,
which is always assumed stationary and ergodic.
We also assume here that the environment is {\it elliptic}, that
is there exists an $\epsilon>0$ such that
$P$-a.s.,
$\omega(x,x+e)\geq \epsilon$ for all $x,e\in  {\ZZ}^d,
|e|=1$. Finally, we denote by $\PP$ the {\it annealed} law
of the RWRE started at $0$,
 that is the law of $\{X_n\}$ under the measure
$P\times  P^0_\om$.

The RWRE model has a natural physical motivation and interpretation
in terms of transport in random media. Mathematically, and especially
for $d>1$,
it leads to
the analysis of irreverservible, inhomogeneous Markov chains, to which
standard tools of homogenization theory do not apply well.
Further, unusual phenomena, such as sub-diffusive
behavior, polynomial decay
of probabilities of large deviations, and trapping effects, arise,
already
in the one dimensional model.

When $d=1$, we write $\omega_x=\omega(x,x+1)$,
$\rho_x=\omega_x/(1-\omega_x)$, and $u=E_P\log \rho_0$.
The following
reveals some of the surprising phenomena associated with the RWRE:
\begin{theorem}[Transience, recurrence, limit speed, {\boldmath $d=1$}]
\label{theo-1} (a) With $\mbox{\rm sign}(0)$ $=1$, it holds that $\PP$-a.s.,
$$\limsup_{n\to\infty}
\kappa X_n=\mbox{\rm sign}(\kappa u) \infty\,,\quad \kappa=\pm 1\,.$$
Further, there is a $v$ such that
\begin{equation}
\label{speed}
\lim_{n\to\infty}\frac{X_n}{n}=v\,, \quad \PP-{\mbox {\rm a.s.}}\,,
\end{equation}
$v>0$ if $\sum_{i=1}^\infty E_P(\prod_{j=0}^i \rho_{-j})<\infty$,
$v<0$ if $\sum_{i=1}^\infty E_P(\prod_{j=0}^i \rho_{-j}^{-1})<\infty$,
and $v=0$ if both these conditions do not hold.\\
(b) If $P$ is a product measure then
\begin{equation}
\label{speed1}
v=\left\{\begin{array}{ll}
\;\;\;\frac{1-E_P(\rho_0)}
{1+E_P(\rho_0)}\,, & E_P(\rho_0)<1,\\
-\,\frac{1-E_P(\rho_0^{-1})}
{1+E_P(\rho_0^{-1})}\,, &
E_P(\rho_0^{-1})<1,\\
\,\;\;\;0\,, & \mbox{\rm else}.
\end{array}
\right.
\end{equation}
\end{theorem}

Theorem \ref{theo-1} is essentially due to \cite{solomon}, see
\cite{stflour} for a proof in the general ergodic setup.
The surprising features of the RWRE model alluded to above can be
appreciated if one notes, already for a product measure $P$, that
the RWRE can be transient with zero speed $v$. Further, if
$P$ is a product measure and $v_0(\omega)$
denotes the speed of a (biased) simple random walk with
probability of jump to the right equal, at any site, to $\omega_0$,
then Jensen's inequality reveals that $|v|\leq |E_P (v_0(\omega))|$,
with examples of strict inequality readily available.

The reason for this behavior is that the RWRE spends a large time in
small traps. This is very well understood in the case $d=1$,
to which the next
section is devoted. We introduce there certain hitting times,
show how they yield precise information on the RWRE, and describe
the analysis of these hitting times.
Understanding the behavior of the RWRE when
$d>1$ is a major challenging problem, on which much progress has
been done in recent years, but for  which many embarrassing
open questions remain. We give a glimpse of what is involved
in Section \ref{section 3},
where  we introduce certain {\it regeneration} times,
and show their usefulness in a variety of situations.
Here is a particularly simple setup where
law of large numbers (and CLT's, although we do not emphasize
that here) are available:
\begin{theorem}
\label{theo-2}
Assume $P$ is a product measure,
$d\geq 6$, and $\omega(x,x+e)=\eta>0$ for $e=
\pm e_i, i=1,\ldots,5$. Then there exists a deterministic constant
$v$ such that $X_n/n\to v$, $\PP$-a.s..
\end{theorem}

\section{The one-dimensional case} \label{section 2}
\setzero\vskip-5mm \hspace{5mm }

\noindent {\bf Recursions}

Let us begin with a sketch of the proof of
Theorem \ref{theo-1}. The transience and recurrence
criterion is proved by noting that conditioned on the environment
$\omega$, the Markov chain $X_n$ is reversible. More explicitly,
fix an interval $[-m_-,m_+]$ encircling the origin
and for $z$ in that interval, define
$$
\calV_{m_-, m_+, \om} (z) :=
P_\om^z (\{X_n\} \text{\ hits\ $-m_-$ before hitting $m_+$)}
\,.
$$
Then,
\begin{equation}
\label{2}
\calV_{m_-, m_+, \om} (z) =
\frac{\dsum_{i=z+1}^{m_+}\;
\dprod_{j=z+1}^{i-1} \rho_j}
{\dsum_{i=z+1}^{m_+} \; \dprod_{j=z+1}^{i-1} \rho_j +
 \dsum_{i=-m_-+1}^z
\left(\dprod_{j=i}^z \rho_j^{-1}\right)}\,,
\end{equation}
from which the conclusion follows.
The proof of the LLN
is more instructive:
define the hitting times
$T_n=\min\{t>0: X_t=T_n\}$, and set $\tau_i=T_{i+1}-T_i$. Suppose
that $\limsup_{n\to\infty} X_n/n=\infty$.  One
checks that $\tau_i$ is an ergodic sequence, hence
$T_n/n\to \EE(\tau_0)$ $\PP$-a.s., which in turns implies that
$X_n/n\to 1/\EE(\tau_0)$, $\PP$-a.s..
But,
$$\tau_0= \won_{\{X_1=1\}}+
\won_{\{X_1=-1\}}(1+\tau_{-1}'+\tau_0')\,,$$
where $\tau_{-1}'$ $(\tau_0')$  denote the first hitting time
of  $0$ ($1$) for the random walk $X_n$ after it  hits $-1$.
Hence, taking $P^0_\om$ expectations,  and noting that
$\{E_{P^i_\om}(\tau_{i})\}_i$ are, $P$-a.s., either all finite or all
infinite,
\begin{equation}
\label{pp1}
  E_{P^0_\om}(\tau_0)=\frac{1}{\omega_0}+
\rho_0  E_{P^{-1}_\om}(\tau_{-1})\,.
\end{equation}
When $P$ is a product measure,
$\rho_0$ and
$E_{P^{-1}_\om}(\tau_{-1})$
are $P$-independent, and taking expectations results with
$\EE(\tau_0)=(1+E_P(\rho_0))/(1-E_P(\rho_0))$
if the right hand side is positive and
$\infty$ otherwise, from which (\ref{speed1}) follows.
The ergodic case is obtained by iterating the relation (\ref{pp1}).

The hitting times $T_n$ are also the beginning of the study of limit
laws for $X_n$. To appreciate this in the case of product
measures $P$ with $E_P(\log \rho_0)<0$
(i.e., when the RWRE is transient to
$+\infty$),  one first observes that
from the above recursions,
$$\EE(\tau_0^r)<\infty \Longleftrightarrow E_P(\rho_0^r)<1 \,.$$
Defining $s=\max\{r: E_P(\rho_0^r)<1\}$, one then expects
that $(X_n-v n)$, suitably rescaled, possesses a limit law,
with $s$-dependent
scaling. This is indeed the case:
 for $s>2$, it is not hard to check that one obtains a
central limit
theorem
with scaling $\sqrt{n}$ (this holds true in fact for ergodic
environments under appropriate mixing
assumptions and with a suitable
definition of the parameter $s$, see \cite{stflour}).
For $s\in (0,1)\cup (1,2)$, one obtains in the i.i.d. environment
case
a Stable($s$) limit law with scaling $ n^{1/s}$ (the cases
$s=1$ or $s=2$ can also be handled but involve logarithmic factors
in the scaling and the deterministic shift). In particular, for
$s<2$ the walk is {\it sub-diffusive}.
We omit the details, referring to \cite{KKS} for
the proof,  except to say that the extension to ergodic
environments of many of
these results has recently been carried out,
see \cite{alex}.

\noindent {\bf Traps}

 The unusual behavior of one dimensional RWRE
is due to the existence of traps in the environment. This is
exhibited most dramatically when one tries to evaluate the
probability of slowdown of the RWRE. Assume that $P$ is a product
measure, $X_n$ is transient to $+\infty$ with positive speed $v$
(this means that $s>1$ by Theorem \ref{theo-1}), and that
$s<\infty$ (which means that $P(\omega_0<1/2)>0$). One then has:
\begin{theorem}[\cite{DPZ,GZ}]
\label{theo-3}
For any $w\in [0, v)$, $\eta>0$, and  $\delta > 0$ small enough,
\begin{equation}
\limn \frac{\log \PP \left(\frac{X_n}{n} \in (w-\delta, w +
\delta)\right)}{\log n} = 1-s\,,
\end{equation}
\begin{equation}
\label{subexpof1} \liminf_{n\to\infty} \frac{1}{n^{1-1/s+\eta}}
\log P^0 \left(\frac{X_n}{n} \in (w-\delta, w + \delta)\right) = 0
\,,\quad P-a.s.,
\end{equation}
and
\begin{equation}
\label{subexpof3}
\limsup_{n\to\infty} \frac{1}{n^{1-1/s-\eta}}
\log P^0 \left(\frac{X_n}{n}\in (w-\delta,w+\delta)\right)
 = -\infty  \,,\quad P-a.s..
\end{equation}
\end{theorem}
(Extensions of Theorem \ref{theo-3} to the mixing environment
setup are presented in \cite{stflour}. There are also precise
asymptotics available in the case $s=\infty$ and $P(\om_0=1/2)>0$,
see \cite{PP,PPZ}).

One immediately notes the difference in scaling between the
annealed and
quenched slowdown estimates in Theorem \ref{theo-3}.
These are due to the fact that, under the quenched
measure, traps are almost surely of a maximal given
size, determined by $P$, whereas under the annealed measure
$\PP$ one can create, at some cost in probability, larger
traps.

To demonstrate
the role of traps in the RWRE model, let us exhibit, for
$w=0$,  a lower
bound that captures the correct behavior in the  annealed setup, and
that forms the basis for the proof of the more general statement.
Indeed, $\{X_n\leq \delta\}
\subset \{T_{n\delta}\geq n\}$. Fixing $R_k=R_k(\om)
:=k^{-1} \sum_{i=1}^k
\log \rho_i$, it holds that $R_k$ satisfies a large deviation
principle with rate function $J(y)=\sup_{\lambda}
(\lambda y- \log E_P(\rho_0^\lambda))$,
and it is not hard to check that
$s=\min_{y\geq 0} y^{-1} J(y)$. Fixing a $y$ such that
$J(y)/y\leq s+\eta$, and $k=\log n/y$, one checks that
the probability that there exists in $[0,\delta n]$ a point
$z$ with $R_k\circ \theta^z \omega \geq y$ is
at least $n^{1-s-\eta}$. But, the probability that the
RWRE does not cross such a segment by time $n$ is, due to (\ref{2}),
bounded away from $0$ uniformly in $n$. This yields
the claimed lower bound in the annealed case. In the quenched
case, one has to work with traps of size almost
$k=\log n/ sy$ for which $kR_k\geq y$, which occur
with probability 1 eventually, and use (\ref{2}) to compute
the probability of an atypical slowdown inside such a trap.
The fluctuations in the length of these typical traps is the
reason why the slowdown probability is believed, for $P$-a.e.
$\omega$, to fluctuate with $n$, in the sense that
$$\liminf_{n\to\infty} \frac{1}{n^{1-1/s}}
\log P^0_\om \left(\frac{X_n}{n} \in (-\delta,
\delta)\right) = -\infty \,,\quad P-a.s.,
$$
while it is known that
$$\limsup_{n\to\infty} \frac{1}{n^{1-1/s}}
\log P^0_\om \left(\frac{X_n}{n} \in (-\delta,
\delta)\right) = 0 \,,\quad P-a.s..
$$
This has been demonstrated rigorously in some
particular cases, see \cite{Gnode}.

The role of traps, and the difference  they produce
between the quenched and annealed regimes, is dramatic also in
the scale of large deviations. Roughly, the
exponential (in $n$) rate of decay of the probability
of atypical events differ between  the quenched and
annealed regime:
\begin{theorem}
\label{theo-4}
The random variables $X_n/n$ satisfy, for
$P$-a.e. realization of the environment $\omega$,
a large deviations principle (LDP) under $P_\om^0$
with a deterministic rate function $I_P(\cdot)$.
Under the annealed measure $\PP$, they satisfy a LDP
with rate function
\begin{equation}
\label{trap}
I(w) = \inf_{Q\in {\cal M}_1^e} (h(Q|P)+I_Q(w))\,,
\end{equation}
where $h(Q|P)$ is the  specific entropy of $Q$ with respect
to $P$ and ${\cal M}_1^e$ denotes the space of stationary ergodic
measures on $\Omega$.
\end{theorem}

Theorem \ref{theo-4} means that to create an annealed large
deviation, one may first ``modify'' the environment (at a certain
exponential cost) and then apply  the quenched LDP in the new
environment. We refer to \cite{greven} (quenched) and
\cite{CGZ,DGZ} for proofs and generalizations to non i.i.d.
environments. We also note that Theorem \ref{theo-4} stands in
sharp contrast to what happens for
random walks on Galton-Watson trees,
where the growth of the tree creates enough
variability in the (quenched) environment to make the annealed and quenched
LDP's identical, see \cite{DGPZ}.

\noindent {\bf Sinai's recurrent walk and aging}

 When $E_P(\log
\rho_0)=0$, traps stop being local, and the whole environment
becomes a diffused trap. The walk spends most of its time ``at the
bottom of the trap'', and as time evolves it is harder and harder
for the RWRE to move. This is the phenomenum of {\it aging},
captured in the following theorem:
\begin{theorem}
\label{theo-5}
There exists a random variable $B^n$,
depending on the environment only, such
that
$$
\PP \left( \left| \frac{X_n}{(\log n)^2} - B^n \right| > \eta
\right) \underset{n\to\infty}{\to}0
\,.$$
Further, for $h>1$,
\begin{equation}
\label{sinai10}
\lim_{\eta\to0} \lim_{n\to\infty} \PP
\left(\frac{|X_{n^h} - X_n|}{(\log n)^2} < \eta\right)
= \frac{1}{h^2} \left[ \frac{5}{3} - \frac{2}{3} e^{-(h-1)} \right]
\,.\end{equation}
\end{theorem}

The first part of Theorem \ref{theo-5} is due to Sinai
\cite{sinai}, with Kesten \cite{kesten} providing the evaluation
of the limiting law of $B^n$. The second part is implicit in
\cite{golosov2}, we refer to \cite{alice} and \cite{stflour} for
the proof and references.

\section{Multi-dimensional RWRE}
 \label{section 3} \setzero\vskip-5mm \hspace{5mm }

\noindent {\bf Homogenization}

 Two special features simplify the
analysis of the RWRE in the one-dimensional case: first, for every
realization of the environment, the RWRE is a reversible Markov
chain. This gave transience and recurrence criteria. Then, the
location of the walk at the hitting times $T_n$ is deterministic,
leading to  stationarity and mixing properties of the sequence
$\{\tau_i\}$ and to a relatively  simple analysis of their tail
properties. Both these features are lost for $d>1$.

A (by now standard) approach to homogenization problems
is to consider the {\it environment viewed from the particle}.
More precisely, with $\theta^x$ denoting the $\ZZ^d$ shift by $x$,
the process $\omega_n=\theta^{X_n}\omega$ is a Markov chain
with state-space  $\Omega$. Whenever the invariant measure of
this chain is absolutely continuous with respect to $P$,
law of large numbers and CLT's can be deduced,
see \cite{kozlov}.
For reversible situations, e.g.
in the ``random conductance model'' \cite{masi},
the invariant measure of the chain $\{\omega_n\}$
is known explicitly. In the non-reversible RWRE model,
this approach has had limited consequences:  one needs to
establish absolute continuity of the invariant measure without
knowing it  explicitly. This was done in \cite{lawler} for
balanced environments, i.e. whenever $\omega(x,x+e)=\omega(x,x-e)$
$P$-a.s. for all $e\in \ZZ^d, |e|=1$, by developing a-priori
estimates on the invariant measure., valid for {\it every} realization
of the environment. Apart from that (and the very recent
\cite{firas}), this approach has not been very useful in the
study of RWRE's.

\noindent {\bf  Regeneration}

 We focus here on another approach
based on analogs of hitting times. Throughout, fix a direction
$\ell\in\ZZ^d$, and consider the process $Z_n=X_n\cdot \ell$.
Define the events $A_{\pm \ell}=\{Z_n\to_{n\to\infty} \pm
\infty\}$. Then, with $P$ a product measure, one shows that
$\PP(A_\ell\cup A_{-\ell})\in \{0,1\}$, \cite{kalikow}. We sketch
a proof:
Call a time $t$ {\it fresh} if $Z_t>Z_n, \forall n<t$, and
for any fresh time $t$,
define the return time $D_t=\min\{n>t: Z_n<Z_t\}$,
calling
$t$ a regeneration time if $D_t=\infty$.
Then, $\PP(A_\ell)>0$ implies by the Markov property
that $\PP(A_\ell\cap \{D_0=\infty\})>0$.
Similarly, on $A_\ell$, each fresh time has a bounded
away from zero probability
to be a regeneration time. One deduces
that  $\PP(\exists\,\mbox{\rm a regeneration time}|A_\ell)=1$.
In particular, on $A_{\pm \ell}$, $Z_n$ changes signs only finitely
many times. If $\PP(A_\ell\cup A_{-\ell})<1$ then with positive
probability, $Z_n$ visits a finite centered interval
infinitely often, and
hence it must
change signs infinitely many times.
But this implies that $\PP(A_\ell\cup A_{-\ell})=0$.

The proof above can be extended to non-product $P$-s having good mixing
properties using, due to the uniform ellipticity, a coupling with
simple nearest neighbor random walk. This is done as follows:
 Set $W = \{0\}\cup
\{\pm e_i\}_{i=1}^d$.
Define the measure
$$
\oo\PP = P \otimes Q_\eps \otimes \oo{P}_{\om,{\beps}}^0
\quad \mbox{\rm
on}\quad
\Bigl( \Om \times W^\NN \times (\ZZ^d)^\NN\Bigr)
$$
in the following way:
$Q_\eps$ is a product measure, such that with $\beps
=(\eps_1,\eps_2,\ldots)$
denoting an element of $W^\NN$,
$Q_\eps (\eps_1 = \pm e_i) = \eps/2$, $i=1, \cdots, d$,
$Q_\eps(\eps_1=0) = 1- \eps d$.
For each fixed $\om,\beps$, $\oo{P}_{\om,\beps}^0$ is the
law of the Markov chain $\{X_n\}$ with state space $\ZZ^d$, such that
$X_0=0$ and, for each $e\in W$, $e\not= 0$,
$$
\oo{P}_{\om,\beps}^0 (X_{n+1} = z+e | X_n = z) =
\won_{\{\eps_{n+1} =e\}} + \frac{\won_{\{\eps_{n+1}=0\}}}{1-d\eps}
[\om (z,z+e)-\eps/2]\,.
$$
It is not hard to check that the law of $\{X_n\}$ under
$\oo{\PP}$ coincides with its law under $\PP$, while its law under
$Q_\eps \otimes \oo{P}_{\om,\beps}^0$ coincides
with its law under $P^0_\om$.
Now, one introduces modified regeneration times
$D_t^{(L)}$ by requiring that after the fresh time $t$,
the ``$\beps$'' coin was used for $L$ steps in the direction $\ell$:
more precisely, requiring that
$\eps_{t+i}=u_i, i=1, \ldots, L$ for some fixed sequence
$u_i\in \ZZ^d, |u_i|=1, u_i\cdot\ell>0$ such that $\sum_{i=1}^L u_i\cdot
\ell\geq L/2$. This, for large $L$, introduces enough
decoupling to carry through the proof,
see \cite[Section 3.1]{stflour}.
We can now state the:\\
{\bf Embarrassing Problem 1}\, Prove that $\PP(A_\ell)\in
\{0,1\}$.

For $d=2$, and $P$ i.i.d., this was shown in \cite{zernermerkl},
where counter examples using non uniformly elliptic, ergodic $P$'s
are also provided. The case $d>2$, even for $P$ i.i.d., remains
open.

\noindent {\bf Embarrassing Problem 2}\, Find transience and
recurrence criteria for the RWRE under $\PP$.

The most promising approach so far
toward Problem 2 uses
regeneration times. Write $0\leq d_1<d_2<\ldots$ for the ordered
sequence of
regeneration times, assuming that $\PP(A_\ell)=1$.
The  name regeneration time is justified
by the following property, which for simplicity we state
in the case $\ell=e_1$:
\begin{theorem}[\cite{sznitmanzerner}]
\label{reg}
For $P$ a product measure,
the sequence \\
$$\{\{\om_z\}_{z\cdot\ell\in [Z_{d_i},Z_{(d_{i+1}-1)})},
\{X_t\}_{t\in [d_i,d_{i+1})}\}_{i=2, 3,\ldots}$$ is i.i.d..
\end{theorem}

From this statement, it is then not hard to deduce that
once $\EE(d_2-d_1)<\infty$, a law of large numbers results, with
a non-zero limiting velocity.
Sufficient conditions for transience put forward in
\cite{kalikow} turn out to fall in this class,
see \cite{sznitmanzerner}. More recently, Sznitman
has introduced a condition that ensures both a LLN and a CLT:

\noindent
{\it Sznitman's T' condition:
$\PP(A_\ell)=1$ and , for some $c>0$ and all $
\gamma<1$, $$\EE(\exp (c \sup_{0\leq n<d_1} |X_n|^\gamma))<
\infty.$$}

\noindent
A remarkable fact about Sznitman's T' condition is that he was able
to derive, using renormalization techniques,
a (rather complicated) criterion, depending on the restriction of $P$
to finite boxes, to check it. Further, Sznitman's T' condition
implies a
good control on $d_1$, and in particular that $d_1$ possesses
all moments, which is the key to the LLN and CLT statements:
$$\EE \left(\exp \left(\log d_1\right)^\delta\right)<\infty,
\forall \delta<2d/(d+1)\,.$$
For these, and related, facts see \cite{sznitmancriterion}.
This leads one to the

\noindent {\bf Challenging Problem 3}\, Do there exist
non-ballistic RWRE's for $d>1$ satisfying that $\PP(A_\ell)=1$ for
some $\ell$?

For $d=1$, the answer is affirmative, as we saw, as soon as $ E_P
\log \rho_0<0$ but $s<1$. For $d>1$, one suspects that the answer
is negative, and in fact one may suspect that $\PP(A_\ell)=1$
implies Sznitman's condition T'. The reason for the striking
difference is that for $d>1$, it is much harder to force the walk
to visit large traps.

It is worthwhile to note that the modified regeneration times
$\{D_t^{(L)}\}$
can be used
to deduce the LLN for a class of mixing environments. We refer to
\cite{CZ} for details. At present, the question of CLT's in such
a general set up remains open.

\noindent {\bf Cut points}

Regeneration times are less useful if the walk is not ballistic.
Special cases of non-ballistic models have been analyzed in the
above mentioned \cite{lawler}, and using a heavy renormalization
analysis, in \cite{bric} for the case of symmetric, low disorder,
i.i.d. $P$. In both cases, LLN's with zero speed  and CLT's are
provided. We now introduce, for another special class of models, a
different class of times that are not regeneration times but
provide enough decoupling to lead to useful consequences.

The setup is similar to that
in Theorem \ref{theo-2}, that is we assume that $d\geq 6$ and that
the RWRE, in its first 5 coordinate, performs a deterministic
random walk:\\
$ \mbox{\rm For $i=1,\ldots,5$}\,, \quad
\om(x,x\pm e_i)=q_{\pm i}\,, \quad
\mbox{\rm for some deterministic
$q_{\pm i}$}\,, \quad P-a.s..
$\\
Set $S=
\sum_{i=1}^5 (q_i+q_{-i})$,
let $\{R_n\}_{n\in \ZZ}$
denote a (biased) simple random walk in $\ZZ^5$ with
transition probabilities $q_{\pm i}/S$, and fix
a sequence
of independent Bernoulli random variable with $P(I_0=1)=
S$,  letting $U_n=\sum_{i=0}^{n-1} I_i$.
Denote by $X_n^1$ the first 5 components of $X_n$ and by
$X_n^2$ the remaining components.
Then, for every realization $\om$,
the RWRE $X_n$ can be constructed
as the Markov chain with
$X_n^1=R_{U_n}$ and
transition
probabilities
$$\oo P^0_\om(X^2_{n+1}=z|X_n)=\left\{
\begin{array}{ll}
1, & X_n^2=z, I_n=1,\\
\om(X_n,(X_n^1,z))/(1-S), & I_n=0\,.
\end{array}
\right.
$$
Introduce now, for the walk $R_n$, cut times
$c_i$ as those times where the past and future of the path $R_n$ do
not intersect. More precisely,
with ${\cal P}_{I}=\{X_n\}_{n\in I}$,
$$c_1=\min\{t\geq 0: {\cal P}_{(-\infty,t)}\cap {\cal P}_{[t,\infty)}
=\emptyset\}\,,
c_{i+1}=
\min\{t> c_i: {\cal P}_{(-\infty,t)}\cap {\cal P}_{[t,\infty)}
=\emptyset\}\,.$$
The cut-points sequence depends on the ordinary random
walk $R_n$ only. In particular, because that walk evolves in $\ZZ^5$,
it follows, as in
\cite{erdos},
 that there are infinitely many cut points,
and moreover that they have a positive density.
Further,
the increments
$X_{c_{i+1}}^2-X_{c_i}^2$ depend on disjoint parts of the environment.
Therefore, conditioned on $\{R_n,I_n\}$, they are independent
if $P$ is a product measure, and they possess good mixing properties
if $P$ has good mixing properties.  From here, the statement
of Theorem \ref{theo-2} is not too far. We refer the reader
to \cite{BSZ}, where this and CLT
statements (with 5  replaced by a larger integer) are proved. An amusing
consequence of \cite{BSZ} is that for $d>5$, one may
construct ballistic RWRE's with, in the notations
of Section \ref{section 2}, $E_P(v_0(\omega))=0$!

\noindent {\bf Challenging Problem 4}\, Construct cut points for
``true'' non-ballistic RWRE's.

 The challenge here is to construct
cut points and prove that their density is positive, without
imposing a-priori that certain components of the walk evolve
independently of the environment.

\noindent {\bf Large deviations}

We conclude the discussion of multi-dimensional RWRE's by
mentioning large deviations for this model. Call a RWRE {\it
nestling} if $\mbox{\rm co}\,\mbox{\rm supp} Q$, where $Q$ denotes
the law of $\sum_{e\in \ZZ^d: |e|=1} e \om(0,e)$. In words, an
RWRE is nestling if by combining local drifts one can arrange for
zero drift. One has then:
\begin{theorem}[\cite{zerner}]
\label{theo-zerner}
Assume $P$ is a product nestling measure. Then, for $P$-almost every
$\omega$, $X_n/n$ satisfies a LDP under $P_\om^0$ with
deterministic rate function.
\end{theorem}

The proof of Theorem \ref{theo-zerner} involves hitting times: let
$T_{y}$ denote the first hitting time of $y\in \ZZ^d$. One then
checks, using the subaddititve ergodic theorem, that
$$
\Lambda(y,\lambda):=
\lim_{n\to\infty} n^{-1}\log
E_\om^0(\exp (-\lambda T_{ny})\won_{\{T_{ny}<\infty\}})$$
exists and is deterministic, for $\lambda\geq 0$. In the
nestling regime,
where slowdown has sub-exponential decay rate due to the
existence of traps much as for $d=1$, this and concentration
of measure estimates are  enough to yield
the LDP. But:

\noindent {\bf Embarrassing Problem 5}\, Prove the quenched LDP
for non-nestling RWRE's.

 A priori, non nestling walks should have been
easier to handle than nestling walks due to good control on the
tail of regeneration times!

\noindent {\bf Challenging Problem 6}\, Derive an annealed LDP
for the RWRE, and relate the rate function to the quenched one.

 One does not expect a relation as simple as in Theorem \ref{theo-4},
because the RWRE can avoid traps by contouring them, and to change
the environment in a way that surely modifies the behavior of the
walk by time $n$ has probability which seems to decay at an
exponential rate faster than $n$. This puts the muti-dimensional
RWRE in an intermediate position between the one-dimensional RWRE
and walks on Galton-Watson trees \cite{DGPZ}. We also note that
certain estimates on large deviations for RWRE's, without matching
constants, appear in \cite{sznitman-slow}.

\label{lastpage}

\end{document}